\documentclass[preprint]{imsart}

\RequirePackage{amsthm,amsmath,amsfonts,amssymb}
\RequirePackage[numbers]{natbib}

\RequirePackage[OT1]{fontenc}
\RequirePackage[numbers]{natbib}
\RequirePackage[colorlinks,citecolor=blue,urlcolor=blue]{hyperref}
\usepackage{caption}
\usepackage{dsfont}
\usepackage{subcaption}
\usepackage{epstopdf}
\usepackage{bm}

\startlocaldefs

\startlocaldefs
\numberwithin{equation}{section}

\newtheorem{theorem}{Theorem}[section]

\newtheorem{lemma}[theorem]{Lemma}

\newcommand{\beq}{\begin{equation}}  
\newcommand{\eeq}{\end{equation}}

\newcommand{\ben}{\begin{enumerate}}  
\newcommand{\een}{\end{enumerate}}
\newcommand{\bed}{\begin{itemize}}  
\newcommand{\eed}{\end{itemize}}

\newcommand {\ee} {{\mathbb E}}
\newcommand{\E}{\mathbb{E}}
\newcommand {\FF} {{\mathcal F}}
\newcommand {\GG} {{\mathcal G}}

\newcommand{\rr}{\mathbb{R}}

\newcommand {\TT} {{\mathcal T}}

\def\qt#1{\qquad\text{#1}}

\newcommand{\one}{\mathds{1}}

\endlocaldefs

\endlocaldefs

\begin{document}

\begin{frontmatter}
\title{Minimax bounds for estimating multivariate Gaussian
  location mixtures}
\runtitle{Minimax bounds for multivariate Gaussian location mixtures}

\begin{aug}
\author[A]{\fnms{Arlene K. H.} \snm{Kim}\ead[label=e1]{arlenent@korea.ac.kr}},
\and
\author[B]{\fnms{Adityanand} \snm{Guntuboyina}\ead[label=e2]{aditya@stat.berkeley.edu}}

\runauthor{Kim and Guntuboyina}

\address[A]{Department of Statistics, 
Korea University, 
145 Anam-ro, Seongbuk-gu,
Seoul, 02841, South Korea \\
 \printead{e1}}

\address[B]{Department of Statistics,
423 Evans Hall,
University of California,
Berkeley, CA - 94720, USA \\
 \printead{e2}}
\end{aug}

\begin{abstract}
We prove minimax bounds for estimating Gaussian location
  mixtures on $\mathbb{R}^d$ under the squared $L^2$ and the
  squared Hellinger loss functions. Under the squared $L^2$ loss, we
  prove that the minimax rate is upper and lower bounded by a constant
  multiple of $n^{-1}(\log n)^{d/2}$. Under the squared Hellinger
  loss, we consider two subclasses based on the behavior of the tails
  of the mixing measure. When the mixing measure has a sub-Gaussian
  tail, the minimax rate under the squared Hellinger loss is bounded
  from below by $(\log n)^{d}/n$. On the other hand, when the mixing
  measure is only assumed to have a bounded $p^{\text{th}}$ moment for
  a fixed $p > 0$, the minimax rate under the squared Hellinger loss
  is bounded from below by $n^{-p/(p+d)}(\log n)^{-3d/2}$. These rates
  are minimax optimal up to logarithmic factors. 
\end{abstract}


\end{frontmatter}


\section{Introduction}
Let $\phi$ be the standard univariate normal density and, for $d \geq
1$, let $\FF_d$ denote the class of densities on $\rr^d$ of the form:   
\begin{align}
(x_1, \ldots, x_d) &\mapsto \int \phi(x_1-u_1) \phi(x_2-u_2) \ldots
                     \phi(x_d-u_d) dG(u_1, \ldots, u_d), \label{class} 
\end{align}
where $G$ is a probability measure on $\rr^d$. $\FF_d$ is precisely
the class of all Gaussian location mixture densities on $\rr^d$. We
study minimax rates in the problem of estimating an unknown density
$f^* \in \FF_d$ from i.i.d observations $X_1, \ldots, X_n$ (throughout
the paper, we assume that $n \geq 2$). 

The minimax rate crucially depends on the choice of the loss
function. We study two different loss functions in this paper. The
first is the squared $L^2$ distance:
\begin{equation}\label{l2dist}
L^2(f, g) :=  \int \left(f(\bm x) - g(\bm x) \right)^2 d\bm x. 
\end{equation}
The minimax risk of estimation over $\FF_d$ under the $L^2$ loss
function is
\begin{equation*}
  R_n\left(\FF_d, L^2 \right) := \inf_{\hat{f}_n} \sup_{f^* \in
 \FF_d} \E_{f^*} L^2 \left(\hat{f}_n, f^* \right). 
\end{equation*}
In Theorem \ref{l2rate}, we prove that $R_n\left(\FF_d, L^2 \right)$
is of the order $n^{-1}(\log n)^{d/2}$. This result is known for $d =
1$. Indeed, when $d = 1$, the upper bound follows from the 
results proved in \citet{ibragimov2001} for estimation of smooth
functions (also see \citet[Section 3]{kim2014}) and the lower bound
was proved by \citet{kim2014}. To the best of our knowledge the result 
for $d \ge 2$ is novel. It is interesting that the rate $n^{-1}(\log
n)^{d/2}$ has a relatively mild dependence on the dimension $d$ and
thus the usual curse of dimensionality is largely avoided for
estimating multivariate Gaussian location mixtures under the $L^2$
loss function.   

The second loss function we investigate is the squared Hellinger
distance: 
\begin{equation}\label{h2dist}
h^2(f,g) := \int \left(\sqrt{f(\bm x)}-\sqrt{g(\bm x)} \right)^2 d \bm x.
\end{equation}
In order to obtain meaningful rates under the
squared Hellinger distance, it is necessary to impose additional
conditions on the probability measure $G$ underlying the density
\eqref{class}. The most common assumption in the literature is to
assume that $G$ is discrete with a known upper bound on the number of
atoms. The Hellinger accuracy (as well as accuracy in the total
variation distance) of estimating discrete Gaussian location mixtures
has been investigated, for example, in \cite{DWYZ2020, HN2016,
  LS2017, SOAJ2014, WY2020}. In particular, it was proved by
\citet{DWYZ2020} (and \citet{WY2020} for $d = 1$) that the minimax
rate is $n^{-1}$ when the dimension $d$ and the number of atoms of $G$
are bounded from above by constants.

In contrast to the discrete mixture situation, minimax rates in
squared Hellinger distance under broader assumptions on $G$ are
not fully understood. Given a subclass $\GG$ of probability measures
on $\rr^d$, let $\FF_\GG$ denote the class of all densities of the
form \eqref{class} where $G$ is constrained to be in $\GG$. We shall 
denote the minimax risk over $\FF_{\GG}$ in the squared Hellinger
distance by 
\[
R_n(\FF_\GG, h^2) := \inf_{\hat f_n} \sup_{f^* \in \FF_\GG} \E_{f^*} h^2(\hat f_n, f^*).
\] We study $R_n(\FF_\GG, h^2)$ for the following two natural subclasses $\GG$:  
\ben
\item $\GG = \GG_1(\Gamma)$: the class of all probability measures $G$
  satisfying $G\{\bm u: \|\bm u\|>t \} \leq \Gamma \exp(-t^2/\Gamma)$ for all $t
  > 0$ and a constant $\Gamma$. Every probability measure in
  $\GG_1(\Gamma)$ has sub-Gaussian tails.    
\item $\GG = \GG_2(p,K)$: the class of all probability measures $G$ satisfying 
\[
\left(\int  \|\bm u\|^p dG(\bm u) \right)^{1/p} \leq K
\]
for a fixed $p>0$ and constant $K > 0$. 
\een
The problem of estimation of densities belonging to the classes
$\FF_{\GG_1(\Gamma)}$ and $\FF_{\GG_2(p, K)}$ has been studied in
\cite{GvdV2001, Zhang2009, kim2014} for $d = 1$ and in \cite{saha2017}
for $d \geq 1$. Extending the results of \citet{Zhang2009}  to $d \geq
1$, \citet{saha2017} analyzed the performance of the nonparametric
maximum likelihood estimator over $\FF_d$ leading to the following 
upper bounds on $R_n(\FF_{\GG_1(\Gamma)}, h^2)$ and $R_n(\FF_{\GG_2(p,
  K)}, h^2)$: 
\begin{equation}\label{ub.sg}
R_n(\FF_{\GG_1(\Gamma)}, h^2) \leq  C_{d, \Gamma} \frac{(\log n)^{d+1}}{n} ,
\end{equation}
and
\begin{equation}\label{ub.mom}
R_n(\FF_{\GG_2(p, K)}, h^2) \leq C_{d,K,p} n^{-\frac{p}{p+d}} (\log n)^{\frac{2d+2p+dp}{2p+2d}}.
\end{equation}
To the best of our knowledge, the corresponding lower bounds do not
currently exist in the literature (except for the case of 
$R_n(\FF_{\GG_1(\Gamma)}, h^2)$ for $d = 1$) and we establish these in
this paper. Specifically, we prove that
\begin{equation}\label{lb.sg}
  R_n(\FF_{\GG_1(\Gamma)}, h^2) \geq  c_{d, \Gamma} \frac{(\log n)^{d}}{n} 
\end{equation}
and
\begin{equation}\label{lb.mom}
  R_n(\FF_{\GG_2(p, K)}, h^2) \geq c_{d,K,p} n^{-\frac{p}{p+d}} (\log n)^{-\frac{3d}{2}}
\end{equation}
in Theorem \ref{lb_subgaussian} and Theorem \ref{lb_pth}
respectively. 

\eqref{lb.sg} implies that there is a logarithmic price to be paid
for dimensionality under the sub-Gaussianity
assumption. \eqref{lb.mom} implies that the rate of convergence
becomes much slower (than the parametric rate) if we only assume
boundedness of the $p^{th}$ moment of $G$ for a fixed $p > 0$. It is
usually believed that Gaussian location mixtures are arbitrarily
smooth leading to nearly parametric rates of 
estimation. While this is true for the $L^2$ loss function, our
results reveal that the story is more complicated for the squared
Hellinger loss function. Specifically, inequality \eqref{lb.mom} shows
that, under the squared Hellinger distance, the rates can be
arbitrarily slow if the mixing measure is allowed to have heavy
tails. This fact does not seem to have been emphasized previously in
the literature even for $d = 1$. Furthermore, for each fixed $p$, the rate
becomes exponentially slow in $d$ revealing the usual curse of
dimensionality.  Note that our lower bounds also imply that the
upper bounds \eqref{ub.sg} and \eqref{ub.mom} cannot be substantially
improved. 

The rest of the paper is organized as follows. Our main results are
all stated in the next section. Theorem \ref{l2rate}  proves the
minimax rate of $(\log n)^{d/2}/n$ for $\FF_d$ under the $L^2$
loss. Theorem \ref{lb_subgaussian} and Theorem \ref{lb_pth} deal with
the squared Hellinger loss function: Theorem \ref{lb_subgaussian}
proves the minimax lower bound of $(\log n)^d/n$ under the 
subgaussianity assumption on the mixing measure and Theorem
\ref{lb_pth} proves the $n^{-p/(p+d)} (\log  n)^{-3d/2}$ lower bound
under the bounded $p^{th}$ moment assumption on the mixing
measure. The proofs of these results are given in Section
\ref{sec:proofs}. We also recall in this section (see Subsection
\ref{prelims}) some basic facts about Fourier transforms, Hermite
polynomials and Assouad's lemma that are used in our proofs.

\section{Main Results}
We state all our main results in this section. Our first result shows
that $R_n(\FF_d, L^2)$  is of the order $n^{-1} (\log n)^{d/2}$. 

\begin{theorem}\label{l2rate}
  There exist constants $c_d$ and $C_d$ depending only on $d$ such that
  \begin{equation}\label{l2rate.eq}
c_d\frac{(\log n)^{d/2}}{n} \leq    R_n(\FF_d, L^2) \leq C_d
\frac{(\log n)^{d/2}}{n}.  
  \end{equation}
\end{theorem}

The proof of the upper bound on $R_n(\FF_d, L^2)$ in \eqref{l2rate.eq}
is based on an extension of the ideas of \cite{ibragimov2001} to $d
\ge 1$ (a simple exposition of these ideas can be found in
\citet[Theorem 4.1]{kim2014}). It involves considering the estimator
\begin{equation}\label{sinc.est}
  \hat{f}_n(\bm x) := \frac{1}{nh^d} \sum_{i=1}^n \mathbb{K} \left(\frac{X_i
  - \bm x}{h}\right) 
\end{equation}
where $\mathbb{K}(\bm y) := K(y_1) \ldots K(y_d)$ with $K(y) :=
(\sin y)/(\pi y)$ and the bandwidth $h$ is taken to be $h := (2\log
n)^{-1/2}$. Controlling the variance of $\hat f_{n}(\bm
  x)$ is straightforward while bounding the bias is non-trivial and we
  do this via Fourier analysis.  

The proof of the lower bound on $R_n(\FF_d, L^2)$ in \eqref{l2rate.eq}
is based on an extension of the ideas of \citet[Proof of Theorem
1.1]{kim2014}. It involves 
applying Assouad's Lemma (recalled in Lemma \ref{assouad}) to a
carefully chosen subset of $\mathcal{F}_d$ whose elements
are indexed by a hypercube. This subset of $\mathcal{F}_d$ is
constructed by taking mixing measures that are additive perturbations
of a Gaussian mixing measure. The additive perturbations are created
using Hermite polynomials.

Our next result proves a lower bound of order $(\log n)^d/n$ for
$R_n(\FF_{\mathcal{G}_1(\Gamma)}, h^2)$. A comparison with the upper
bound \eqref{ub.sg} of \citet{saha2017} reveals that this lower bound
is possibly off by at most a factor of $\log n$ and is thus minimax
rate optimal up to the single $\log n$ multiplicative factor.

\begin{theorem}\label{lb_subgaussian}
There exist a positive constant $c_{d}$ depending only on $d$ such
that 
\[
R_n (\FF_{\GG_1(\Gamma)}, h^2) \geq c_{d} \frac{(\log n)^d}{n} \qt{for
all $\Gamma \geq c^{-1}_d$}. 
\]
\end{theorem}


The proof of Theorem 2.2 is based on an extension
  of the ideas of \citet[Proof of Theorem 1.3]{kim2014}. Assouad's
  lemma is applied to a subset of $\FF_{\GG_1(\Gamma)}$ which is constructed by
  taking mixing measures that are additive perturbations of a Gaussian
  mixing measure. The perturbations are different from those used in
  the proof of Theorem \ref{l2rate} although they are also based on
  Hermite polynomials. It is not easy to 
  directly work with the squared Hellinger loss function while dealing
  with mixture densities so we crucially use the fact that the squared
  Hellinger loss is bounded from below by a constant multiple of the
  chi-squared divergence for the constructed subset of $\FF_d$.

Our final result proves a lower bound of the order $n^{-p/(p+d)}$
(up to a logarithmic factor) for $R_n(\FF_{\GG_2(p, K)}, h^2)$. As we
mentioned previously, this result reveals that rates strictly slower
than $n^{-1}$ are possible for estimating Gaussian location mixtures
and that the rate can also be affected by the usual curse of
dimensionality.  A comparison with the upper bound \eqref{ub.mom} of
\citet{saha2017} reveals that this lower bound is optimal up to
logarithmic factors possibly depending on $d$. 

\begin{theorem}\label{lb_pth}
For every $p > 0$, there exists a positive constant $c_{d, p}$
depending only on $d$ and $p$ such that  
\[
R_n(\FF_{\GG_2(p, K)}, h^2) \geq c_{d, p}  n^{-\frac{p}{p+d}} (\log
n)^{-\frac{3d}{2}} \qt{for all $K \geq c_{d, p}^{-1}$}. 
\]
\end{theorem}


For the proof of theorem \ref{lb_pth}, we first construct a normal
mixture density whose mixing measure is a discrete
distribution that is supported on a $d$-dimensional product set
(lattice) and has a Pareto type tail behavior. We then construct a
hypercube of normal mixture densities by perturbing this support
set. For each point in the support set, we use either the original
point or a nearby point whose distance to the non-perturbed
  point is determined by the probability value at the original support
  point of the mixing distribution. We finally apply Assouad's
lemma to the constructed hypercube of densities in $\FF_{\GG_2(p,
  K)}$. 

The proofs of our results are given in the next section. 

\section{Proofs}\label{sec:proofs}

\subsection{Preliminaries}\label{prelims}
We shall recall here some standard facts about the Fourier transform
and Hermite polynomials that we shall use in our main proofs. 

We use the notation $\TT$ for the Fourier transform. It is defined as 
\[
\mathcal{T}(f) (\bm{t}) =(2 \pi)^{-d/2}\int
e^{-it_1x_1 - i t_2 x_2 - \dots - i t_d x_d} f( \bm{x}) d \bm{x}
\]
for functions $f \in L_1(\mathbb R^d)$. The Fourier inversion theorem
gives
\begin{equation*}
  f(x) = (2 \pi)^{-d/2} \int e^{i x_1 t_1 + \dots + i x_d t_d}
  \TT(f)(\bm t) d\bm t.
\end{equation*}
Plancherel's theorem states that  
\[
\int | f( \bm{x}) |^2 d \bm{x} = \int | \mathcal{T}(f) ( \bm{t}) |^2 d \bm{t}
\]
for functions $f \in L_1(\mathbb R^d) \cap L_2(\mathbb R^d)$.
The convolution-product property of the Fourier transform states that
\begin{equation*}
  \TT(h)(\bm t) = \TT(f)(\bm t) \TT(g)(\bm t) ~ \text{ for all } \bm t
  ~ \text{where } h(\bm x) := \int f(\bm x - \bm u) g(\bm u) d\bm u.  
\end{equation*}

Our upper bound on $R_n(\FF_d, L^2)$  is based on the sinc kernel
$\mathbb{K}$ with $\mathbb{K}(\bm y) := K(y_1) \ldots K(y_d)$ where
\begin{equation*}
  K(y) := \frac{\sin y}{\pi y}. 
\end{equation*}
It is well-known that
\begin{equation}\label{sinc.fourier}
  K(y) = K(-y), ~~~~ \TT(K)(t) = \frac{1}{\sqrt{2 \pi}} \one_{\{|t| \leq
  1\}}, ~~~~ \text{ and } ~~~~ \TT(K^2)(t) = \frac{1}{\pi \sqrt{2 \pi}}
\left(1 - \frac{|t|}{2} \right)_+
\end{equation}
where $x_+ := \max(x, 0)$. 

Our lower bound constructions involve Hermite polynomials. These are
defined for $d = 1$ as
\begin{equation}\label{hermite}
  H_j(x) := (-1)^j e^{x^2/2} \frac{d^j}{d x^j} e^{-x^2/2} =
  \frac{(-1)^j}{\phi(x)} \frac{d^j \phi(x)}{d x^j}. 
\end{equation}
Note that $H_j(\cdot)$ is an odd function when $j$ is odd (and even
when $j$ is even). The Hermite polynomials are orthogonal with respect
to the weight function $\phi(x) := (2 \pi)^{-1/2}
\exp(-x^2/2)$. Specifically for all $j \geq 0$ and $k \geq 0$, we have
\begin{equation}\label{hermite.orthogonal}
\int \phi(x) H_j(x) H_{k}(x) dx = 
  \begin{cases} 
   0 & \text{for } j \neq k \\
   j!       & \text{for } j = k.
  \end{cases}
\end{equation}

We shall use the bound
\begin{equation}\label{hermite.bound}
  |H_j(x)| \leq \kappa \sqrt{j!} \exp(x^2/4)
\end{equation}
for a constant $\kappa \sim 1.086 < 2^{1/4}$ (see e.g. Equation 8.954
of \cite{table}). 

For $d \geq 1$, we take
\[
H_{ \bm{j}} ( \bm{x}) := \frac{(-1)^{| \bm{j}|}}{\bm\phi( \bm{x})} \left(
  \frac{\partial}{\partial  \bm{x}} \right)^{ \bm{j}} \bm\phi( \bm{x}) =
H_{j_1}(x_1) \ldots H_{j_d}(x_d) 
\]
where $\bm\phi(\bm{x}) = \phi(x_1) \ldots \phi(x_d)$, $\left(
  \partial/\partial  \bm{x} \right)^{ \bm{j}} :=  \partial^{j_1 +
  \ldots + j_d} / \partial x_1^{j_1} \ldots \partial x_d^{j_d}$ and
$|\bm{ j}| = j_1+\ldots + j_d$. 

We next recall Assouad's lemma (see e.g., \cite[Lemma
24.3]{vaart}) which will be our main tool for proving
minimax lower bounds. 
\begin{lemma}[Assouad]\label{assouad}
  Let $d^2$ be either \eqref{l2dist} or \eqref{h2dist} and define the
  minimax risk
  \begin{equation*}
    R_n(\FF_d, d^2) := \inf_{\hat f_n} \sup_{f^* \in \FF_d} \E_{f^*}
    d^2 \left(\hat f_n, f^* \right). 
  \end{equation*}
For some $N \geq 1$, let $\{f_{\bm \tau}, \bm\tau \in \{0, 1\}^N\}$ be a
subset of $\FF_d$. Then
\begin{equation}\label{assouad.eq}
  R_n(\FF_d, d^2) \geq \frac{N}{8} \min_{\bm\tau \neq \bm\tau'}
  \frac{d^2(f_{\bm\tau}, f_{\bm\tau'})}{\Upsilon(\bm\tau, \bm\tau')}
  \min_{\Upsilon(\bm\tau, \bm\tau') = 1} \left(1 - \sqrt{\frac{n}{2}
      \chi^2(f_{\bm\tau} \| f_{\bm\tau'})} \right)
\end{equation}
where $\Upsilon(\bm\tau, \bm\tau') := \sum_{j=1}^N \one_{\{\tau_j \neq
  \tau'_j\}}$ is the Hamming distance between $\bm\tau$ and $\bm\tau'$ and
$\chi^2(f \| g) := \int (f - g)^2/g$ is the $\chi^2$-divergence between
densities $f$ and $g$.  
\end{lemma}

\subsection{Proof of Theorem \ref{l2rate}}

We break this proof into two parts: the upper bound on $R_n(\FF_d,
L^2)$ and the lower bound on $R_n(\FF_d, L^2)$. Let us first prove the
upper bound.

We consider the kernel estimator \eqref{sinc.est} with bandwidth $h =
(2 \log n)^{-1/2}$. We need to control its variance and squared
bias. We first bound the variance as
\begin{align}
\ee_{f} \int (\hat f_n(\bm y) - \ee \hat f_n(\bm y))^2 d\bm y &=
                                                                \frac{1}{n
                                                                h^{2d}}
  \int \text{var}_f \left(\mathbb{K} \left(\frac{X_1 - y}{h} \right)
                                                                \right)d\bm
                                                                y \nonumber \\ &\leq
                                                                 \frac{1}{nh^{2d}}
                                                                 \int
                                                                 \ee_{f}
                                                                 \mathbb{K}^2
                                                                 \left(\frac{X_1
                                                                 - y}{h} \right)
                                                                d\bm
                                                                 y
                                                                 \nonumber\\  
&= \frac{1}{nh^{2d}} \int \int \prod_{i=1}^d
                                                                                 K^2\left(\frac{x_{i}-y_i}{h} \right) f(\bm x) d\bm x d\bm y \nonumber\\
&=   \frac{1}{nh^{2d}} \int \left(\prod_{i=1}^d \int K^2
                                                                                                                                                           \left(\frac{x_i
                                                                                                                                                           -
                                                                                                                                                           y_i}{h}
                                                                                                                                                           \right)
                                                                                                                                                           dy_i\right)
                                                                                                                                                           f(\bm
                                                                                                                                                           x)
                                                                                                                                                           d\bm
                                                                                                                                                           x
  \nonumber\\
                                                              &=
                                                                \frac{1}{nh^{d}}
                                                                \int
                                                                \left(\prod_{i=1}^d
                                                                \int
                                                                K^2(z_i)
                                                                dz_i\right)
                                                                f(\bm
                                                                x)
                                                                d\bm x \nonumber
  \\
  &= \frac{\pi^{-d}}{n h^d} \leq \frac{1}{nh^d}
                                                                \label{vari}
\end{align}
where we used $\int K^2(z) dz = \int |\TT(K)(z)|^2 dz = 1/\pi$ and
$\int f(\bm x) d \bm x = 1$. 

For the bias, note that $\bm y \mapsto \ee_{f^*} \hat
f_n(\bm y)$ has Fourier 
transform  $\bm t \mapsto (2\pi)^{d/2} \TT(f)(\bm t) \prod_{i=1}^d
\TT(K)(ht_i)$. This is because  
\[
\ee_{f} \hat f_n(\bm y) = \frac{1}{h^d} \int \prod_{i=1}^d K \left(\frac{x_i-y_i}{h} \right)f(\bm x) d\bm x 
\]
and
\begin{align*}
\TT(\ee_f \hat f_n) (\bm t) &= \frac{1}{(2\pi)^{d}h^d} \int e^{-i {\bm t}^T \bm y} \left[ \prod_{i=1}^d K \left(\frac{x_i-y_i}{h} \right)f(\bm x) d\bm x\right] d\bm y \\
&= \frac{1}{(2\pi)^{d/2}h^d} \int \prod_{i=1}^d \left( \int e^{-it_i y_i} K \left(-\frac{x_i-y_i}{h} \right) dy_i\right) f(\bm x) d\bm x \\
&=  \frac{1}{(2\pi)^{d/2}h^d} \int \prod_{i=1}^d \left( \int e^{-it_i (y_i-x_i)} K \left(\frac{y_i-x_i}{h} \right) dy_i\right)e^{-i\bm t^T \bm x} f(\bm x) d\bm x \\
&= \prod_{i=1}^d \left(\int \frac{1}{\sqrt{2\pi}} e^{-i(ht_i)z_i} K(z_i) dz_i \right)  \int  e^{-i\bm t^T \bm x} f(\bm x) d\bm x \\
&= \prod_{i=1}^d \TT(K)(ht_i) (2\pi)^{d/2} \TT(f)(\bm t).
\end{align*}
Thus by Plancherel theorem (and \eqref{sinc.fourier})
\begin{align*}
\text{bias}^2 = \int (\ee_f \hat f_n(\bm y) - f(\bm y))^2 &= \int \left| \TT(\ee_f \hat f_n)(\bm t) - \TT(f)(\bm t) \right|^2 d\bm t \\
&= \int |\TT(f)(\bm t)|^2 \left| (2\pi)^{d/2} \prod_{i=1}^d \TT(K)(ht_i) -1 \right|^2 d\bm t\\
&= \int |\TT(f) (\bm t)|^2 |\one_{\{|t_1|\leq 1/h, \ldots, |t_d| \leq 1/h  \}}-1|^2 d\bm t \\
&\leq \sum_{i=1}^d \int_{|t_i| > 1/h}  |\TT(f) (\bm t)|^2 d\bm t,
\end{align*}
where the last inequality follows since $|\one_{\{|t_1|\leq 1/h,
  \ldots, |t_d| \leq 1/h  \}}-1|^2 \leq \sum_{i=1}^d \one_{\{|t_i| >
  1/h\}}$. The density $f \in \FF_d$ is of the form \eqref{class} for
some probability measure $G$ and thus we can assume that $f$ is the
density of $X = U + Z$ for independent random variables $U \sim G$
and $Z$ having the standard normal distribution on $\mathbb{R}^d$. We
can thus write
\begin{align*}
  \TT(f)(\bm t) &= (2 \pi)^{-d/2} \E \exp(- i \left<\bm t, U + Z
  \right>) \\ &= (2 \pi)^{-d/2} \E \exp(-i \left<\bm t, U \right>) \E
  \exp(-i \left<\bm t, Z \right>)  
\end{align*}
Using the formula for the characteristic function of the Gaussian random
variable $Z$ and the fact that the characteristic function of $U$ is
bounded by 1, we get
\begin{equation*}
  |\TT(f)(\bm t)| \leq (2 \pi)^{-d/2} \exp \left(-\frac{1}{2}
    \sum_{i=1}^d t_i^2 
  \right). 
\end{equation*}
As a result,
\begin{equation*}
 \text{bias}^2 \leq \sum_{i=1}^d \int_{|t_i| > 1/h} |\TT(f)(\bm
 t)|^2 d \bm t \leq (2 \pi)^{-d} \sum_{i=1}^d \int_{|t_i| > 1/h}  \exp \left(-\sum_{i=1}^d t_i^2 
  \right) d\bm t. 
\end{equation*}
The standard Gaussian tail bound $\int_{|u| > a} e^{-u^2} du \leq 2
\sqrt{\pi} e^{-a^2}$ now leads to
\begin{equation*}
   \text{bias}^2 \leq \pi^{-d/2} 2^{-d-1} d \exp\left(\frac{-1}{h^2}
   \right). 
\end{equation*}
Combining this with \eqref{vari}, we get 
\begin{equation*}
  R_n(\FF_d, L^2) \leq \inf_{h > 0} \left(\frac{1}{nh^d} + \pi^{-d/2} 2^{-d-1} d \exp\left(\frac{-1}{h^2}
   \right) \right). 
\end{equation*}
The choice $h := (\log n)^{-1/2}$ then clearly leads to $R_n(\FF_d,
L^2) \leq C_d (\log n)^{d/2}/n$ which proves the upper bound. 

We now prove the lower bound on $R_n(\FF_d,
L^2)$. The idea is to construct a subset of $\FF_d$ indexed by a
hypercube $\{0, 1\}^N$ for some $N$ and then use Assouad's lemma
(Lemma \ref{assouad}). Our construction is a natural extension to $d
\geq 1$ of the one-dimensional construction in \cite{kim2014} and is
described below. Let $m$ be the largest integer such that
\begin{equation}\label{mcond}
  \frac{m^{5d/4} 8^d 3^{dm}}{\sqrt{n}} \le 1.
\end{equation}
We can assume without loss of generality that $n$ is large enough so
that $m$ defined as above satisfies $m \geq 3$. It is also easy to 
check (because $m^{5d/4} 8^d 3^{dm} \geq e^m$) that the above
condition for $m$ implies that 
\begin{equation*}
  m \leq \frac{1}{2} \log n. 
\end{equation*}
Below we denote by
$\phi_{\sigma^2}(\cdot)$, the univariate normal density with mean zero
and variance $\sigma^2$ i.e., $\phi_{\sigma^2}(x) := (\sqrt{2 \pi}
\sigma)^{-1} \exp(-x^2/(2 \sigma^2))$. We construct densities (for the
application of Assouad's lemma) via perturbations of the density:
\begin{equation*}
  f_0(\bm x) := \int \phi(x_1 - u_1) \dots \phi(x_d - u_d) \gamma(\bm u)
  d\bm u \qt{where $\gamma(\bm u) := \phi_m(u_1) \dots \phi_m(u_d)$}. 
\end{equation*}
Note also that
\begin{equation*}
  f_0(\bm x) = \phi_{1+m}(x_1) \dots \phi_{1+m}(x_d). 
\end{equation*}
Now let
\begin{equation*}
  \mathcal{J} := \{1, 3, \dots, 2m-1\}^d
\end{equation*}
and note that cardinality of $\mathcal{J}$ is $|\mathcal{J}| =
m^d$. We shall apply Assouad's lemma (Lemma \ref{assouad}) with $N :=
m^d$ and we index binary vectors in $\{0, 1\}^N$ by elements $\bm j :=
(j_1, \dots, 
j_d)$ of $\mathcal{J}$. For $\bm \alpha = (\alpha_{\bm j}, \bm j \in
\mathcal{J})$, let
\begin{equation*}
  f_{\bm \alpha}(\bm x) = \int \phi(x_1 - u_1) \dots \phi(x_d - u_d)
  \left[\gamma(\bm u) + \epsilon \sum_{\bm j \in \mathcal{J}}
    \alpha_{\bm j} \gamma_{\bm j}(\bm u)  \right] d\bm u 
\end{equation*}
where $\epsilon$ is given by
\begin{align}\label{eps.def}
  \epsilon = c n^{-1/2} m^{-d/4} 
\end{align}
for a constant $c \in (0, 1)$ that will be determined later and  
\begin{equation*}
  \gamma_{\bm j}(\bm u) := \gamma_{j_1}(u_1) \dots \gamma_{j_d}(u_d)
  \qt{with } \gamma_{j_i}(u_i) := 2^{1/2}(2\pi)^{3/4}
  \sqrt{\frac{3^{j_i}}{j_i!}} \phi(u_i)
  H_{j_i}\left(\frac{2}{\sqrt{3}}u_i \right). 
\end{equation*}
Here $H_{j_i}(\cdot)$ denotes the Hermite polynomial (see
\eqref{hermite}). Let us first argue that $f_{\bm \alpha} \in
\FF_d$. To see 
this, it is enough to show that
\begin{equation}\label{mix.me}
  \bm u \mapsto \gamma(\bm u) + \epsilon \sum_{\bm j \in \mathcal{J}}
  \alpha_{\bm j} \gamma_{\bm j}(\bm u) 
\end{equation}
integrates to $1$ over $\bm u \in
\mathbb{R}^d$ and is nonnegative. Integration to 1 is justified by the
fact that $\gamma_{j_i}(u_i)$ is an odd function of $u_i$ for each $i = 1,
\dots, d$ and the fact that $\int \gamma(\bm u) d \bm u =
1$. For nonnegativity of \eqref{mix.me}, note first that the
inequality \eqref{hermite.bound} implies that, for each $i = 1, \dots,
d$,  
\begin{align}
  |\gamma_{j_i}(u_i)| &= \sqrt{2} (2 \pi)^{3/4}
                        \sqrt{\frac{3^{j_i}}{j_i !}} \phi(u_i)
                        \left|H_{j_i}(2u_i/\sqrt{3}) \right| \leq \sqrt{2} (2 \pi)^{1/4} 2^{1/4}
    3^{j_i/2} \exp \left(-\frac{u_i^2}{6} \right). \label{boundgamma} 
\end{align}
Because we have assumed that $n$ is large enough so that $m \geq 3$,
we have $\exp(-u_i^2/6) \leq \exp(-u_i^2/(2m))$ which gives
\begin{align*}
  |\gamma_{j_i}(u_i)| \leq 8  \cdot 3^{j_i/2} \sqrt{m} \phi_m(u_i) \qt{for
  every $i = 1, \dots, d$}
\end{align*}
and consequently
\begin{equation*}
|\gamma_{\bm j}(\bm u)| \leq 8^d 3^{|\bm j|/2} m^{d/2} \prod_{i=1}^d
\phi_m(u_i) = 8^d 3^{|\bm j|/2} m^{d/2} \gamma(\bm u)  \qt{where $|\bm
  j| := j_1 + \dots + j_d$}.  
\end{equation*}
Thus
\begin{align*}
  \gamma(\bm u) + \epsilon \sum_{\bm j \in \mathcal{J}} \alpha_{\bm j}
  \gamma_{\bm j} (\bm u) &\geq \gamma(\bm u) \left[1 - \epsilon
                           8^d m^{d/2}\sum_{\bm j \in \mathcal{J}}
                           3^{|\bm j|/2} \right] \\ &\geq \gamma(\bm u) \left[1 - \epsilon
                           8^d m^{d/2} 3^{d(2m-1)/2} |\mathcal{J}| \right]
                                                      = \gamma(\bm u) \left[1 - \epsilon
                           8^d m^{3d/2} 3^{dm}\right]
\end{align*}
because the maximum value of any $j_i$ is $2m-1 \leq 2m$ for every
$\bm j \in \mathcal{J}$ and the cardinality of $\mathcal{J}$ is
$m^d$. Plugging in our value of $\epsilon$ (from \eqref{eps.def}) and
using condition \eqref{mcond}, we get  
\begin{align}\label{lb.mix}
  \gamma(\bm u) + \epsilon \sum_{\bm j \in \mathcal{J}} \alpha_{\bm j}
  \gamma_{\bm j} (\bm u) &\geq \gamma(\bm u) \left( 1 - c n^{-1/2}
                           m^{5d/4} 8^d 3^{dm}\right) \geq \left(1 -
                           c\right) \gamma(\bm u) 
\end{align}
which implies nonnegativity of \eqref{mix.me} as long as $c < 1$.   

We now lower bound $\min_{\bm \alpha \neq \bm \beta} \frac{L^2(f_{\bm
    \alpha}, f_{\bm \beta})}{\Upsilon(\bm \alpha, \bm \beta)}$ (where 
$\Upsilon(\bm \alpha, \bm \beta) := \sum_{\bm j \in \mathcal{J}}
\one_{\{\alpha_{\bm j} \neq \beta_{\bm j}\}}$ is the Hamming distance
between $\bm \alpha$ and $\bm \beta$). Observe first that
\begin{equation*}
  f_{\bm \alpha}(\bm x) - f_{\bm \beta}(\bm x) = \epsilon \int
  \phi(x_1 - u_1) \dots \phi(x_d - u_d) \left\{\sum_{\bm j \in
      \mathcal{J}} \left(\alpha_{\bm j} - \beta_{\bm j} \right)
    \gamma_{\bm j}(\bm u) \right\} d\bm u = \epsilon \sum_{\bm j \in
    \mathcal{J}}\left(\alpha_{\bm j} - \beta_{\bm j} \right)
  \Gamma_{\bm j}(\bm x) 
\end{equation*}
where
\begin{equation*}
  \Gamma_{\bm j}(\bm x) := \int \phi(x_1 - u_1) \dots \phi(x_d - u_d)
  \gamma_{\bm j}(\bm u) d\bm u. 
\end{equation*}
As a result
\begin{align*}
  \int \left(f_{\bm \alpha}(\bm x) - f_{\bm \beta}(\bm x) \right)^2
  d\bm x &= \epsilon^2 \int \left(\sum_{\bm j \in \mathcal{J}}
    (\alpha_{\bm j} - \beta_{\bm j}) \Gamma_{\bm j}(\bm x) \right)^2
    d\bm x\\
     &= \epsilon^2 \int \left(\sum_{\bm j \in \mathcal{J}}
    (\alpha_{\bm j} - \beta_{\bm j}) \TT(\Gamma_{\bm j})(\bm t) \right)^2
    d\bm t. 
\end{align*}
Because $\Gamma_{\bm j}$ is defined as the product over $i = 1, \dots,
d$ of the convolution of $\phi$ and $\gamma_{j_i}$, we have
\begin{equation*}
\mathcal{T} (\Gamma_{\bm j}) (\bm t) = \prod_{i=1}^{d}
2^{1/2}(2\pi)^{3/4}\left[ \sqrt{\frac{3^{j_i}}{j_i!}}
  \mathcal{T}(\phi)(t_i) \mathcal{T}\left(\phi(\cdot)
    H_{j_i}(\frac{2}{\sqrt{3}} \cdot)\right)(t_i)\right].   
\end{equation*}
By \citet[Lemma 2.1]{kim2014}, we have
\begin{equation*}
  \mathcal{T} \left(\phi(\cdot)
    H_{k}(\frac{2}{\sqrt{3}} \cdot)\right)(t) = (-{\rm i})^k 3^{-k/2}
  \phi(t) H_{k}(2t) \qt{for odd $k$}. 
\end{equation*}
This gives (note that $\sqrt{2} (2 \pi)^{3/4} \phi^2(t) = \sqrt{2
  \phi(2 t)}$) 
\begin{equation*}
  \mathcal{T} (\Gamma_{\bm j}) (\bm t)  = \prod_{i=1}^d (- {\rm
    i})^{j_i} \sqrt{2 \phi(2 t_i)} \frac{H_{j_i}(2 t_i)}{\sqrt{j_i!}}
\end{equation*}
and thus
\begin{equation*}
    \int \left(f_{\bm \alpha}(\bm x) - f_{\bm \beta}(\bm x) \right)^2
  d\bm x = \epsilon^2 \int \left(\sum_{\bm j \in \mathcal{J}}
    (\alpha_{\bm j} - \beta_{\bm j}) \prod_{i=1}^d (- {\rm
    i})^{j_i} \sqrt{2 \phi(2 t_i)} \frac{H_{j_i}(2 t_i)}{\sqrt{j_i!}} \right)^2
    d\bm t. 
  \end{equation*}
  The orthogonality of the Hermite polynomials with respect to the
  weight function $\phi$ (see Equation \eqref{hermite.orthogonal}) implies
  \begin{align*}
    \int \left(f_{\bm \alpha}(\bm x) - f_{\bm \beta}(\bm x) \right)^2
  d\bm x &= \epsilon^2 \sum_{\bm j \in \mathcal{J}}  (\alpha_{\bm j} -
  \beta_{\bm j})^2 \int \left(\prod_{i=1}^d (- {\rm 
    i})^{j_i} \sqrt{2 \phi(2 t_i)} \frac{H_{j_i}(2 t_i)}{\sqrt{j_i!}} \right)^2
           d\bm t \\ &= \epsilon^2 \sum_{\bm j \in \mathcal{J}}  (\alpha_{\bm j} -
  \beta_{\bm j})^2 = \epsilon^2 \Upsilon(\bm \alpha, \bm \beta). 
  \end{align*}
  We thus have
  \begin{equation*}
    \min_{\bm \alpha \neq \bm \beta} \frac{L^2(f_{\bm \alpha}, f_{\bm
        \beta})}{\Upsilon(\bm \alpha, \bm \beta)} \geq \epsilon^2. 
  \end{equation*}
Now we bound the $\chi^2$ distance between $f_{\bm \alpha}$ and
$f_{\bm \beta}$ for $\bm \alpha$ and $\bm \beta$ with $\Upsilon(\bm
\alpha, \bm \beta)=1$. Note first that, as a result of \eqref{lb.mix},
we have
\begin{align*}
  f_{\bm \beta}(\bm x) \geq (1 - c) f_0(\bm x)  \qt{for all $\bm \beta$}
\end{align*}
so that
\begin{align*}
  \chi^2\left(f_{\bm \alpha} \| f_{\bm \beta} \right) \leq
  \frac{1}{1-c} \int \frac{\left(f_{\bm \alpha} - f_{\bm \beta}
  \right)^2}{f_0}. 
\end{align*}
We now split the integral above into $R(\bm x):=\{|x_1| \leq Mm^{1/2},
\ldots, |x_d| \leq Mm^{1/2}\} $ and $R(\bm x)^c$ where $M$ is a
dimensional constant larger than $8d \log 3$. Then using
$f_0(\bm x)\{x \in R(\bm x)\} \geq  m^{-d/2}/C_d$ (for some
dimensional constant $C_d$), we
have   
\begin{align*}
(1-c)\chi^2(f_{\bm \alpha}\|f_{\bm \beta}) &\leq \int \frac{(f_{\bm \alpha}(\bm x)-f_{\bm \beta}(\bm x))^2}{f_0(\bm x)} \\
&\leq C_dm^{d/2} \int (f_{\bm \alpha}(\bm x)-f_{\bm \beta}(\bm x))^2+ \int_{R(x)^c} \frac{(f_{\bm \alpha}(\bm x) - f_{\bm \beta}(\bm x))^2}{f_0(\bm x)}d\bm x \\
&= C_d m^{d/2}\epsilon^2 + \int_{R(x)^c} \frac{(f_{\bm \alpha}(\bm
                                                                                                                                                                   x)
                                                                                                                                                                   -
                                                                                                                                                                   f_{\bm
                                                                                                                                                                   \beta}(\bm
                                                                                                                                                                   x))^2}{f_0(\bm
                                                                                                                                                                   x)}d\bm
                                                                                                                                                                   x
                                                                                                                                                                   \leq \frac{1}{2n}   + \int_{R(x)^c} \frac{(f_{\bm \alpha}(\bm
                                                                                                                                                                   x)
                                                                                                                                                                                                                        - f_{\bm \beta}(\bm x))^2}{f_0(\bm x)}d\bm x
\end{align*}
provided $c$ is chosen above so that $c^2 C_d \leq 1/2$ (recall that
$\epsilon^2 = c^2 n^{-1} m^{-d/2}$). The second term above is bounded
as follows. Denoting $\bm j^*$ by the
index where  $\bm \alpha$ and $\bm \beta$ differ so that
$\alpha_{\bm j^*} \neq \beta_{\bm j^*}$ (recall that $\Upsilon(\bm
\alpha, \bm \beta) = 1$), we get by
\eqref{boundgamma}, 
\begin{align*}
\int_{R(x)^c} \frac{(f_{\bm \alpha}(\bm
                                                                                                                                                                   x)
                                                                                                                                                                                                                        - f_{\bm \beta}(\bm x))^2}{f_0(\bm x)}d\bm x &\leq \epsilon^2 3^{d/2} (4\pi)^{3/4} 3^{|\bm j^*|}  \int_{ R(x)^c} \frac{ \Big( \int \phi(x_1-u_1) \ldots \phi(x_d-u_d)  \bm \phi_3(\bm u)  d\bm u\Big)^2} {f_0(\bm x)} d\bm x\\
&\leq d\epsilon^2 3^{d/2} (4\pi)^{3/4} 3^{|\bm j^*|}   \int_{ \{|x_1|> M \sqrt{m}\}} \frac{\bm\phi_{4}^2(\bm x)}{\bm\phi_{1+m}(\bm x)} d \bm x\\
&\leq  2d c^2n^{-1} 3^{d/2} (4\pi)^{3/4} 3^{|\bm j^*|}
                                                                                                                                             \int_{\{x_1>M\sqrt{m}\}}{\phi_{4}(x_1)} dx_1 \leq \frac{1}{2n},
\end{align*}
where the second inequality follows by bounding the maximum by the sum with symmetry,  the penultimate inequality follows first by plugging \eqref{eps.def} in $\epsilon^2$ and then since $\bm\phi_{1+m}(\bm x) \geq m^{-d/2} \bm\phi_{4}(\bm x)$ and $3^{|\bm
  j^*|} \int_{\{ x_1 > Mm^{1/2} \}} \phi_4(x_1) dx_1 \leq 1$ by taking $M$
larger than $8d \log 3$, and the last inequality follows by choosing
$c$ to be a small enough dimensional constant such that $c^2 \leq
(4 \pi)^{-3/4}3^{-d/2}/(2d)$.   

By Assouad's lemma \ref{assouad}, we have 
 \[
R_n(\FF_d, L^2) \geq \frac{m^d}{8} \epsilon^2 \left(
  1-\sqrt{\frac{1}{2(1 - c)}}\right) 
\]
Of course $c$ can be taken to be small enough so that the right hand
side above is larger than a constant (depending on $d$ alone) multiple
of $(\log n)^{d/2}/n$. The proof is thus complete. 

\subsection{Proof of Theorem \ref{lb_subgaussian}}

The idea for the lower bound on $R_n(\FF_{\GG_1(\Gamma)}, h^2)$ is to
again construct a subset of $\FF_{\GG_1(\Gamma)}$ indexed by a
hypercube and then use Assouad's lemma. Our construction is a
natural extension to $d\geq 1$ of the one-dimensional construction in
\cite{kim2014}. Let $m$ be the largest integer such that  
\[
\frac{m^d 6^{d} 5^{dm}}{\sqrt{n}} \leq 1.
\] 
The above condition implies $m \leq \log n$. In order to use Assouad's
lemma, we construct densities via perturbations of the density
\[
f_0(\bm x) = \phi_2(x_1) \ldots \phi_2(x_d)
\]
where, it may be recalled, that $\phi_{\sigma^2}(\cdot)$ denotes the
univariate normal density with mean zero and variance $\sigma^2$. Note
that $f_0$ can also be written as
\[
f_0(\bm x) = \int \phi(x_1-u_1) \ldots \phi(x_d-u_d) v(\bm u) d\bm u \ \ \ \text{where } v(\bm u) := \phi(u_1) \ldots \phi(u_d).
\]
Again we let 
\[
\mathcal{J} := \{1,3, \ldots, 2m-1\}^d
\]
with $N:= m^d$ and we consider index binary vectors in $\{0,1\}^N$ by elements $\bm j:= (j_1, \ldots, j_d)$ of $\mathcal{J}$. For $\bm \alpha = (\alpha_{\bm j}, \bm j \in \mathcal{J})$, we let
\[
f_{\bm{\alpha}}( \bm{x}) = \int \phi(x_1-u_1) \ldots \phi(x_d-u_d) \left[ v(\bm u) + \epsilon \sum_{\bm j \in \mathcal{J}} \alpha_{\bm j} v_{\bm{j}}( \bm{u}) \right] d\bm u
\]
where $\epsilon$ is given by 
\begin{equation}\label{eps2.def}
\epsilon = cn^{-1/2}
\end{equation}
for a constant $c \in (0, 1)$ that will be determined later and  
\[
v_{\bm j} (\bm u) := v_{j_1}(u_1) \ldots v_{j_d} (u_d) \ \ \ \text{with } 
v_{ j_i}( u_i) = 2^{5/4}\sqrt{\pi} \sqrt{\frac{5^{j_i}}{j_i!}} \phi(\sqrt{3}u_i) H_{j_i}\Big( \frac{4}{\sqrt{5}} u_i\Big).
\]
By an analogous argument in the proof of Theorem \ref{l2rate}, we can show $f_{\bm \alpha} \in \FF_{\GG_1(\Gamma)}$. First, integration of $f_{\bm \alpha}(\bm x)$ to 1 is guaranteed since $v_{j_i}(u_i)$ is an odd function of $u_i$ for each $i=1, \ldots, d$ and the fact that $\int v(\bm u) d\bm u=1$. For nonnegativity of $v(\bm u) + \epsilon \sum_{\bm j \in \mathcal{J}} \alpha_{\bm j} v_{\bm j}(\bm u)$, note that using \eqref{hermite.bound} we have
\[
|v_{j_i}(u_i)| = 2^{5/4} \sqrt{\pi} \sqrt{\frac{5^{j_i}}{j_i!}} \phi(\sqrt{3}u_i) \left|H_{j_i}(4u_i/\sqrt{5}) \right| \leq 2 \cdot 5^{j_i/2} \exp \left(-\frac{7u_i^2}{10} \right) \leq 6 \cdot 5^{j_i/2}\phi(u_i)
\]
for every $i=1, \ldots, d$. Thus
\[
|v_{\bm j}(\bm u)| \leq 6^d 5^{|\bm j|/2} v(\bm u)
\]
which gives
\[
v(\bm u) + \epsilon \sum_{\bm j \in \mathcal{J}} \alpha_{\bm j} v_{\bm j}(\bm u) \geq v(\bm u) \left[1-\epsilon 6^d  \sum_{\bm j \in \mathcal{J}} 5^{|\bm j|/2} \right] \geq v(\bm u) \left[ 1-\epsilon 6^d m^{d} 5^{dm} \right].
\]
Plugging in $\epsilon = cn^{-1/2}$ from \eqref{eps2.def}, 
\[
v(\bm u) + \epsilon \sum_{\bm j \in \mathcal{J}} \alpha_{\bm j} v_{\bm j}(\bm u) \geq v(\bm u) \left( 1-cn^{-1/2}6^d m^d 5^{dm} \right) \geq (1-c) v(\bm u) >(1/2) v(\bm u)>0
\]
as long as $c<1/2$. In the same way, we have
\[
v(\bm u) + \epsilon \sum_{\bm j \in \mathcal{J}} \alpha_{\bm j} v_{\bm j}(\bm u) \leq \frac{3}{2} v(\bm u).
\] 
To claim $f_{\alpha} \in \FF_{\GG_1(\Gamma)}$, we need to show $v(\bm u) + \epsilon \sum_{\bm j \in \mathcal{J}} \alpha_{\bm j} v_{\bm j}(\bm u)$ is sub-gaussian.
Indeed, 
\begin{align*}
\int_{\|u\|>t} v(\bm u) + \epsilon \sum_{\bm j \in \mathcal{J}} \alpha_{\bm j} v_{\bm j}(\bm u) d\bm u &\leq \frac{3}{2} \int_{\|u\|>t} \phi(u_1) \ldots \phi(u_d) d\bm u \leq c_d \exp(-t^2/c_d)
\end{align*}
where $c_d$ is a constant depending on $d$.
The fact that all these constructed mixing densities are between $(1/2)v(\bm u)$ and $(3/2) v(\bm u)$ gives 
\begin{equation}\label{f_alpha_0}
\frac{1}{2} f_0(\bm x) \leq f_{\bm \alpha}(\bm x) \leq \frac{3}{2} f_0(\bm x).
\end{equation}
Inequality \eqref{f_alpha_0} implies that 
\begin{equation}\label{hell_chi}
h^2(f_{\bm \alpha}, f_{\bm \beta}) \geq \frac{1}{6} \int \frac{({f_{\bm{\alpha}}}-{f_{\bm{\beta}}})^2}{f_0}
~~ \text{ and } ~~ \chi^2(f_{\bm{\alpha}}\| f_{\bm{\beta}}) \leq  2
\int \frac{({f_{\bm{\alpha}}}-{f_{\bm{\beta}}})^2}{f_0} 
\end{equation}
because, respectively,  
\[
h^2(f_{\bm{\alpha}}, f_{\bm{\beta}}) = \int \frac{({f_{\bm{\alpha}}}-{f_{\bm{\beta}}})^2}{(\sqrt{f_{\bm{\alpha}}}+\sqrt{f_{\bm{\beta}}})^2} \geq \frac{1}{6} \int \frac{({f_{\bm{\alpha}}}-{f_{\bm{\beta}}})^2}{f_0}
\]
and
\[
\chi^2(f_{\bm{\alpha}}\| f_{\bm{\beta}}) = \int \frac{(f_{\bm{\alpha}}-f_{\bm{\beta}})^2}{f_{\bm{\alpha}}} \leq 2 \int \frac{(f_{\bm{\alpha}}-f_{\bm{\beta}})^2}{f_0}.
\]
By inequality \eqref{hell_chi}, it is clear that for the application
of Assouad's Lemma \ref{assouad}, it is enough to focus on the quantity $\int \frac{(f_{\bm{\alpha}}-f_{\bm{\beta}})^2}{f_0}$.  Let us write
\[
\Lambda_{\bm j} := \int \phi(x_1-u_1) \ldots \phi(x_d-u_d) v_{\bm j}(\bm u) d\bm u,
\]
and we consider
\begin{equation}\label{eq1}
\int \frac{(f_{\bm{\alpha}}-f_{\bm{\beta}})^2}{f_0} = \epsilon^2 \int \left(\sum_{ \bm{j} \in \mathcal{J}}(\alpha_{\bm{j}}-\beta_{\bm{j}})  \frac{\Lambda_{\bm j}}{\sqrt{f_0}}\right)^2 = \epsilon^2 \int \left( \sum_{\bm j \in \mathcal{J}} (\alpha_{\bm j}-\beta_{\bm j}) \mathcal{T} \Big(\frac{\Lambda_{\bm j}}{\sqrt{f_0}} \Big) \right)^2 d\bm t. 
\end{equation}
Because $\Lambda_{\bm j}/\sqrt{f_0}$ is defined as the product over $i=1, \ldots,
d,$ of
\[
\frac{\int \phi(x_i-u_i) v_{j_i}(u_i)du_i}{\sqrt{\int \phi(x_i-u_i) \phi(u_i) du_i}} =\sqrt{3\pi}(8\pi)^{1/4}  \sqrt{\frac{ 5^{j_i}}{j_i!}} \int \phi_{\frac{4}{3}}(x_i-u) \phi(\sqrt{\frac{3}{2}}u) H_{j_i}(\frac{3}{\sqrt{5}} u) du,
\]
where the above equality follows by \citet[Lemma 2.2]{kim2014},
we have
\[
\mathcal{T}  \Big(\frac{\Lambda_{\bm j}}{\sqrt{f_0}} \Big)(t_i)  = \prod_{i=1}^d \sqrt{3\pi} (8\pi)^{1/4} \sqrt{\frac{5^{j_i}}{j_i!}} \mathcal{T}(\phi_{4/3})(t_i) \mathcal{T}\left( \phi(\sqrt{\frac{3}{2}} \cdot) H_{j_i}(\frac{3}{\sqrt{5}}\cdot)  \right)(t_i).
\]
By \citet[Lemma 2.1]{kim2014}, we have
\[
\mathcal{T} \left(\phi(\sqrt{\frac{3}{2}} \cdot) H_{j_i} (\frac{3}{\sqrt{5}} \cdot) \right)(t_i) = (-{\rm i})^{j_i} 5^{-j_i/2} \sqrt{\frac{2}{3}} \phi (\sqrt{\frac{2}{3}}t_i) H_{j_i}(2t_i).
\]
Using 
\[
\mathcal{T}(\phi_{4/3})(t_i) = \frac{1}{\sqrt{2\pi}} \exp\left(-\frac{2}{3}t_i^2 \right),
\]
we have 
\[
\mathcal{T}  \Big(\frac{\Lambda_{\bm j}}{\sqrt{f_0}} \Big) = \prod_{i=1}^d (-{\rm i})^{j_i} \sqrt{2\phi(2t_i)} \frac{H_{j_i}(2t_i)}{\sqrt{j_i!}}.
\]
Consequently,
\[
\int \frac{(f_{\bm{\alpha}}-f_{\bm{\beta}})^2}{f_0} = \epsilon^2 \Upsilon(\bm \alpha, \bm \beta).
\]
This (and \eqref{hell_chi}) gives
  \begin{equation*}
    \min_{\bm \alpha \neq \bm \beta} \frac{h^2(f_{\bm \alpha}, f_{\bm
        \beta})}{\Upsilon(\bm \alpha, \bm \beta)} \geq
    \frac{1}{6}\epsilon^2 ~~ \text{ and } ~~   \chi^2(f_{\bm \alpha}, f_{\bm\beta}) \leq 2\epsilon^2.
  \end{equation*}
 By Lemma \ref{assouad}, we have 
 \[
R_n(\FF_{\GG_1(\Gamma)}, h^2)  \geq \frac{m^d}{8} \frac{c^2}{6n} \left(1-c \right) \geq c_d\frac{(\log n)^d}{n}.
\]

%

\subsection{Proof of Theorem \ref{lb_pth}} 
The idea for the lower bound on $R_n(\FF_{\GG_2(p,K)}, h^2)$ is again
to construct a subset of $\FF_{\GG_2(p,K)}$ indexed by a hypercube and
then we use Assouad's lemma. Let $$S := 
\left\{1+(\ell-1)M \log(\ell e),  \ 1 \leq \ell \leq k_0 \right\}$$
where
\begin{equation*}
M  = C_{d,p} \sqrt{\log n} \qt{with $C_{d,p} =
  \sqrt{8\left(\frac{2d}{p+d} +1\right)}$}
\end{equation*}
and $k_0$ the largest integer less than or equal to $k_0'$ where
\begin{equation*}
k_0' = c
  n^{1/(p+d)} (\log 
  n)^{-3/2} \qt{with $c =40^{-1/2}$}. 
\end{equation*}

Let $k=k_0^d$ and let $\bold{a}_1, \dots, 
\bold{a}_k$ be an enumeration of the points in $S^d = S \times \dots
\times S$. We can take $\bold{a}_1=(1, 1,
\ldots, 1) \in \rr^d$, $\bold{a}_2 = (1+M\log(2e),1, \ldots,1),
\ldots, \bold{a}_k = (1+(k_0-1)M\log(k_0e), \ldots,
1+(k_0-1)M\log(k_0e))$. Let $G$ be the discrete probability
distribution supported on $S^d$ that is given by  
\begin{equation}\label{Gdef}
G\{\bold{a}_i\} =C_{n,d,p} \left(\left<\bold{a}_i, \mathbf{1}
  \right> -d+1\right)^{-(p+d)} \qt{for $i = 1, \dots, k$}
\end{equation}
where $\mathbf{1}$ is the $d$-dimensional vector of ones and 
\begin{equation}\label{cdp}
C_{n,d,p} := \frac{1}{\sum_{i=1}^k \left(\left<\bold{a}_i, \mathbf{1}
  \right> -d+1\right)^{-(p+d)} } < \frac{1}{\left(\left<\bold{a}_1, \mathbf{1}
  \right> -d+1\right)^{-(p+d)}} = 1
\end{equation}
is the normalizing constant.
We assume $n$ is sufficiently large so that $\log \log k_0 \leq M$
and $Mp \geq 2^d$. We claim that 
\begin{equation}\label{claim3}
  C_{n,d,p} > 1/2~~ \text{ and } ~~ G \in \GG_2(p,6^{1/p}d). 
\end{equation}
Let us assume the above claim for now and proceed with the proof. The
claim will be proved later. For $i = 1, \dots, k$, let $\bold{b}_i =
\bold{a}_i - (\delta_i/\sqrt{d}) \mathbf{1}$ where 
\[
\delta_i = \frac{1}{\sqrt{2nG(\bold{a}_i)}}.
\]
The choice of $k_0$ means that $\left<\bold{a}_i, \mathbf{1} \right> \leq (d-1)+n^{1/(p+d)}$, which implies that $G(\bold{a}_{i}) > (2n)^{-1}$ for every $i=1, \ldots, k$. Hence $0<\delta_i < 1$.

%
Note that   
$\|\bold{b}_i - \bold{a}_i\|=\delta_i$ for all $i$. 
Define, for every $\bm\alpha = (\alpha_1, \ldots, \alpha_k) \in \{0,
1\}^k$, 
\begin{align*}
f_{\bm \alpha}(\bm{x}) &= \sum_{i=1}^k G(\bold{a}_i)
                         \bm{\phi}_d(\bm{x}-\bold{a}_i(1-\alpha_i)-\bold{b}_i\alpha_i)
\end{align*}
where $\bm{\phi}_d$ is the standard normal density on $\rr^d$.
Since $G \in \GG_2(p,6^{1/p}d)$ and $\|\bold{b}_i\| \leq \|\bold{a}_i\|$ for every $i=1, \ldots, k$,
\[
\sum_{i=1}^k \|(1-\alpha_i) \bold{a}_i + \alpha_i \bold{b}_i \|^p G\{ \bold{a}_i\} \leq  2^{p} \int \|u\|^p dG \leq  6(2d)^p.
\]
This shows $f_{\bm \alpha} \in \FF_{\GG_2(p, 6^{1/p}
  (2d))}$. 

We shall use Assouad's lemma for the class $f_{\bm \alpha}, \bm \alpha
\in \{0, 1\}^k$. For fixed $\bm \alpha = (\alpha_1, \dots, \alpha_k)$
and $\bm \beta = 
(\beta_1, \dots, \beta_k)$ in $\{0, 1\}^k$, we can write 
(with $f_{i}(x) := \bm{\phi}_d(x - a_i(1 -\alpha_i) - b_i \alpha_i)$ and
$g_{j}(x) := \bm{\phi}_d(x - a_j(1 - \beta_j) - b_j \beta_j)$)
\begin{align*}
h^2(f_{\bm\alpha}, f_{\bm\beta}) &= 2 \left( 1-\int \sqrt{f_{\bm\alpha}(\bm{x}) f_{\bm\beta}(\bm{x})} d\bm{x} \right) \\
&= 2\left(1-\int \sqrt{ \sum_{i,j} G(\bold{a}_i)G(\bold{a}_j) f_i(\bm{x})g_j(\bm{x})} d\bm{x} \right)\\
&\geq 2 \left(1-  \sum_{i,j} \sqrt{G(\bold{a}_i)G(\bold{a}_j)}  \int \sqrt{f_i(\bm{x})g_j(\bm{x})}d\bm{x} \right)
\end{align*}
where the last inequality follows from $\sqrt{\sum_i x_i } \leq \sum_i
\sqrt{x_i}$ for $x_i \geq 0$. \citet[Equation (A.18)]{saha2017} now
gives 
\begin{align*}
\frac{1}{2}h^2(f_{\bm\alpha}, f_{\bm\beta}) &=\left( 1- \sum_i
                                              G(\bold{a}_i) \int
                                              \sqrt{f_i(\bm{x})
                                              g_i(\bm{x})} d\bm{x} -
                                              \sum_{i \neq
                                              j}\sqrt{G(\bold{a}_i)G(\bold{a}_j)}\int
                                              \sqrt{f_i(\bm{x})
                                              g_j(\bm{x})}
                                              d\bm{x}\right) \\  
&\geq \sum_{i=1}^{k} G(\bold{a}_i) |\alpha_i-\beta_i| (1-e^{-\delta_i^2/8}) - \sum_{i \neq j} \sqrt{G(\bold{a}_i) G( \bold{a}_j)} e^{-M^2/8}\\
&\geq \sum_{i=1}^k G(\bold{a}_i) \frac{\delta_i^2}{10} |\alpha_i-\beta_i| - \frac{1}{40n},
\end{align*}
where the last inequality follows since $1-e^{-a^2/8} \geq a^2/10$ for every $0<a<1$, and
\[
\sum_{i\neq j} \sqrt{G(\bold{a}_i) G( \bold{a}_j)} e^{-M^2/8} \leq k^2 e^{-M^2/8} = k_0^{2d} e^{-C_{d,p}^2 \log n /8} \leq (1/\sqrt{40})^{2d} n^{\frac{2d}{p+d}-\frac{C_{d,p}^2}{8}} \leq \frac{1}{40n}
\]
by the choice of $k_0$ and $C_{d,p}$. This implies, by the choice of
  $\delta_i^2 = 1/{(2n G(\bold{a}_i))}$, that
\[
h^2(f_{\bm \alpha}, f_{\bm \beta}) \geq  \frac{ \sum_{i=1}^k |\alpha_i-\beta_i| }{20n } - \frac{1}{40n} \geq \frac{\Upsilon(\bm \alpha, \bm \beta)}{40n}.
\]
Now suppose that $\Upsilon(\bm\alpha,\bm\beta)=1$ and let $l$ be
the unique index such that $\alpha_{l} \neq \beta_{l}$. Then for
$\delta_{l}^2 \leq 1$,      
\begin{align*}
\chi^2(f_{\bm \alpha}\| f_{\bm \beta}) &= \int \frac{(f_{\bm \alpha}(\bm x)-f_{\bm \beta}(\bm x))^2}{f_{\bm \alpha}(\bm x)} d\bm x = \int \frac{\left(\sum_{i=1}^k G(\bold{a}_i) f_i(\bm x) - \sum_{i=1}^k G(\bold{a}_i) g_i(\bm x) \right)^2}{\sum_{i=1}^k G(\bold{a}_i) f_i(\bm x)}d\bm x \\
&\leq \int \frac{G(\bold{a}_l)^2 \left( f_l(\bm x) - g_l(\bm x) \right)^2} {G(\bold{a}_l) f_{l}(\bm x)} d\bm x 
= G(\bold{a}_l) \left( \int \frac{g_{l}^2(\bm x)}{f_l(\bm x)}d\bm x -1 \right) \\
&\leq G(\bold{a}_l) \Big(e^{\delta_l^2}-1 \Big)  \leq
                                                                                             2G(\bold{a}_l)\delta_l^2= \frac{1}{n}, 
\end{align*}
where the penultimate inequality follows since 
\[
\int \frac{g_{l}^2(\bm x)}{f_l(\bm x)}d\bm x = \int \frac{\bm{\phi}_d^2(\bm x-\bold{a}_l(1-\beta_l)-\bold{b}_l \beta_l)}{\bm{\phi}_d(\bm x-\bold{a}_l(1-\alpha_l)-\bold{b}_l \alpha_l) } d\bm x \leq \exp \left(\|\bold{a}_l-\bold{b}_l\|^2 \right) = \exp(\delta_l^2),
\] and the last inequality follows since $e^a \leq 2a+1$ when $0\leq a\leq 1$.


By
Lemma \ref{assouad}, we have
\[
R_n(\FF_{\GG_2(p,K)}, h^2) \geq \frac{k}{8} \frac{1}{40n} \left(1-\sqrt{\frac{1}{2}} \right).
\]
Thus the lower bound will be a constant (only depending on $d$ and $p$) multiple of $n^{-p/(p+d)}
(\log n)^{-3d/2}$.

It remains to prove the claim \eqref{claim3}. Let
us first show that $C_{n,d, p} > 1/2$.  For this, let us bound the denominator of $C_{n,d, p}$
using the inequality for $d \geq 2$
\begin{align}
\sum_{\ell_2=1}^{k_0} \cdots \sum_{\ell_d=1}^{k_0} &\Big(1+
\sum_{j=1}^d (\ell_j-1)M \log(\ell_j e)\Big)^{-(p+d)} \leq  \nonumber \\
&(1+(\ell_1-1)M\log(\ell_1 d))^{-(p+1)} \prod_{j=2}^d \left( 1+\frac{1}{M(p+j-1)}
\right). \label{Ad}
\end{align}
The above inequality will be proved at the end of the proof. The above
inequality, along with the inequality
\begin{align*}
\sum_{\ell=1}^{k_0} (1+&(\ell-1)M)^{-(p+1)}
\leq  \left( 1+ \int_{1}^{k_0} (1+(x-1)M )^{-(p+1)} \right) \nonumber \\ 
&=  \left(1+\frac{1}{pM} \Big(1-(1+(k_0-1)M \Big)^{-p} \right) \leq
                                                                            \left(1+\frac{1}{pM} \right), 
\end{align*}
imply that the denominator in the definition \eqref{cdp} of $C_{n,d,p}$
is bounded from above by 
\[
\left(1+\frac{1}{Mp}\right)^d \leq 1+ \frac{(2^d-1)}{Mp} \leq 2.
\]
which proves $1/2 < C_{n,d, p}$. 

We shall next prove that  $\int \| u\|^p dG \leq K$. If $p<1$, then $\int \|u\|^p dG \leq \int \|u\| dG$. Thus without loss of generality, we let $p \geq 1$. By definition of
$G$, 
\begin{align}
\int \| u\|^p dG &=C_{n,d,p} \sum_{u_1, \dots, u_d \in S} \left\{(u_1^2+
                   \ldots + u_d^2)^{p/2} \Big(\sum_{j=1}^d u_j -d+1
                   \Big)^{-(p+d)} \right\} \label{eq2} \\ 
&\leq  d^p \sum_{u_1, \ldots, u_d \in S} \left\{ u_1^p  \Big(\sum_{j=1}^d u_j -d+1 \Big)^{-(p+d)} \right\} \label{eq3}\\
&=  d^p \sum_{u_1 \in S} u_1^p \left\{ \sum_{\ell_2, \ldots, \ell_d} \Big(1+ \sum_{j=1}^d (\ell_j-1)M \log(\ell_j e) \Big)^{-(p+d)} \right\}  \nonumber\\
&\leq  d^p \sum_{u_1 \in S} u_1^{-1} \prod_{j=2}^d
                                                                                                                                  \left(1+\frac{1}{M(p+j-1)} \right), \label{eq4}
\end{align}
where the first inequality follows since $C_{n,d,p}<1$ and $ \sum_{j=1}^d u_j^2  \leq  \left(\sum_{j=1}^d u_j \right)^2 $ for $u_j \geq 1$, which implies $(u_1^2+\ldots+ u_d^2)^{p/2} \leq \left(\sum_{j=1}^d u_j \right)^p$ and this last term is bounded above by $d^{p-1} \sum_{j=1}^d u_j^p$ by H$\ddot{\text{o}}$lder's inequality for $p \geq 1$ so we have \eqref{eq3} by symmetry. Equation \eqref{eq4} follows by the inequality \eqref{Ad}.

It suffices to show $\sum_{u_1 \in S}u_1^{-1}$ is bounded above by a
constant. Indeed,   
\begin{align*}
\sum_{u_1 \in S}u_1^{-1} = \sum_{\ell=1}^{k_0} \left(1+(\ell-1)M \log (\ell e) \right)^{-1} &\leq 1+ \sum_{\ell=2}^{k_0} (1+ (\ell-1) M \log(\ell e))^{-1}\\
& \leq 1+ \sum_{\ell=1}^{k_0-1} (1+\ell M \log (\ell e))^{-1}\\
&\leq 1+ (1+M)^{-1} + \int_1^{k_0} (x M \log x)^{-1} \\
&\leq 2 + \frac{\log \log (k_0)}{M} \leq 3.
\end{align*}
This gives $\int \|u\|^p dG \leq 3d^p \left(1+\frac{1}{Mp}\right)^{d-1} \leq 6d^p$.

The only remaining thing is to prove \eqref{Ad}. For
convenience, for $j' \in \{1, \ldots, d\}$, we let $A_{j'} =
A_{j'}(\ell_1, \ldots, \ell_{j'}) := 1+ \sum_{j=1}^{j'} (\ell_j -1)
M \log(\ell_j e)$. Note that $A_{j'} \geq 1$ for all $j' = 1, \ldots, d$. We need to
show that  for $j' \geq 2$,
\begin{equation}\label{trivialeq}
\sum_{\ell_2}^{k_0}  \cdots \sum_{\ell_{j'}}^{k_0} A_{j'}^{-(p+j')} \leq  A_1^{-(p+1)} \prod_{j=2}^{j'} \left(1+\frac{1}{M(p+j-1)} \right).
\end{equation}
For the above, we use the idea of mathematical induction on $j'$. When $j'=2$, 
\begin{align*}
\sum_{\ell_2=1}^{k_0} A_2^{-(p+2)} &= \sum_{\ell_2=1}^{k_0} \left( 1+(\ell_1-1)M\log(\ell_1 e)+(\ell_2-1)M\log(\ell_2 e) \right)^{-(p+2)} \\
&\leq \sum_{\ell_2=1}^{k_0} \left(A_1+(\ell_2-1)M \right)^{-(p+2)} \\
&\leq A_1^{-(p+2)} + \int_1^{k_0} (A_1 +(x-1)M)^{-(p+2)} dx \\
&\leq A_1^{-(p+1)} \left(1+\frac{1}{M(p+1)} \right)
\end{align*}
where the last inequality follows since $A_1 \geq 1$.
Let the above claim \eqref{trivialeq} is true for $j'=j_0 (\geq 2)$. Then, we consider the case $j'=j_0+1$.
\begin{align*}
\sum_{\ell_2=1}^{k_0}  \cdots &\sum_{\ell_{j_0}=1}^{k_0} \sum_{\ell_{j_0+1}=1}^{k_0} A_{j_0+1}^{-(p+j_0+1)} \leq \sum_{\ell_2=1}^{k_0}  \cdots \sum_{\ell_{j_0}=1}^{k_0} \sum_{\ell_{j'}=1}^{k_0} \left( A_{j_0} + (\ell_{j'}-1)M \right)^{-(p+j_0+1)} \\
&= \sum_{\ell_2=1}^{k_0}  \cdots \sum_{\ell_{j_0}=1}^{k_0}  \left(A_{j_0}^{-(p+j_0+1)} + \int_1^{k_0} \big(A_{j_0} +(x-1)M \big)^{-(p+j_0+1)} dx \right) \\
&\leq \sum_{\ell_2=1}^{k_0}  \cdots \sum_{\ell_{j_0}=1}^{k_0} \left\{ A_{j_0}^{-(p+j_0)}  \left(1+\frac{1}{M(p+j_0)} \right) \right\}\\
&\leq A_1^{-(p+1)} \prod_{j=2}^{j_0} \left(1+\frac{1}{M(p+j-1)} \right) \left(1+\frac{1}{M(p+j_0)} \right)\\
&=A_1^{-(p+1)}  \prod_{j=2}^{j_0+1} \left(1+\frac{1}{M(p+j-1)} \right), 
\end{align*}
where the last inequality follows by the induction hypothesis.

%
%

\section*{Acknowledgements}
%
The first author was supported by NRF-2020R1F1A1A01069632.
The second author was supported by NSF CAREER Grant DMS-16-54589.
%





\end{document}